\documentclass[reqno]{amsart}
\usepackage{amssymb}
\usepackage[pagebackref,hypertex]{hyperref}
\allowdisplaybreaks[4]
\theoremstyle{plain}
\newtheorem{thm}{Theorem}[section]
\newtheorem{lem}{Lemma}[section]
\newtheorem{cor}{Corollary}[section]
\theoremstyle{remark}

\DeclareMathOperator{\td}{d\mspace{-2mu}}
\date{Drafted on 26 November 2008 and completed on 8 December 2008 in VU Student Village}
\date{Revised on 24 November 2010 in Tianjin}
\date{}

\begin{document}

\title{An extension of an inequality for ratios of gamma functions}

\author[F. Qi]{Feng Qi}
\address[F. Qi]{Department of Mathematics, College of Science, Tianjin Polytechnic University, Tianjin City, 300160, China}
\email{\href{mailto: F. Qi <qifeng618@gmail.com>}{qifeng618@gmail.com}, \href{mailto: F. Qi <qifeng618@hotmail.com>}{qifeng618@hotmail.com}, \href{mailto: F. Qi <qifeng618@qq.com>}{qifeng618@qq.com}}
\urladdr{\url{http://qifeng618.wordpress.com}}

\author[B.-N. Guo]{Bai-Ni Guo}
\address[B.-N. Guo]{School of Mathematics and Informatics, Henan Polytechnic University, Jiaozuo City, Henan Province, 454010, China}
\email{\href{mailto: B.-N. Guo <bai.ni.guo@gmail.com>}{bai.ni.guo@gmail.com}, \href{mailto: B.-N. Guo <bai.ni.guo@hotmail.com>}{bai.ni.guo@hotmail.com}}

\begin{abstract}
In this paper, we prove that for $x+y>0$ and $y+1>0$ the inequality
\begin{equation*}
\frac{[\Gamma(x+y+1)/\Gamma(y+1)]^{1/x}}{[\Gamma(x+y+2)/\Gamma(y+1)]^{1/(x+1)}} <\biggl(\frac{x+y}{x+y+1}\biggr)^{1/2}
\end{equation*}
is valid if $x>1$ and reversed if $x<1$ and that the power $\frac12$ is the best possible, where $\Gamma(x)$ is the Euler gamma function. This extends the result in [Y. Yu, \textit{An inequality for ratios of gamma functions}, J. Math. Anal. Appl. \textbf{352} (2009), no.~2, 967\nobreakdash--970.] and resolves an open problem posed in [B.-N. Guo and F. Qi, \emph{Inequalities and monotonicity for the ratio of gamma functions}, Taiwanese J. Math. \textbf{7} (2003), no.~2, 239\nobreakdash--247.].
\end{abstract}

\keywords{extension, inequality, ratio, gamma function}

\subjclass[2010]{Primary 33B15; Secondary 26D07}

\thanks{The first author was partially supported by the China Scholarship Council and the Science Foundation of Tianjin Polytechnic University}

\thanks{This paper was typeset using \AmS-\LaTeX}

\maketitle

\section{Introduction}

It is common knowledge that the classical Euler gamma function $\Gamma(x)$ may be defined for a real argument $x>0$ by
\begin{equation}\label{egamma}
\Gamma(x)=\int^\infty_0t^{z-1} e^{-t}\td t.
\end{equation}
The logarithmic derivative of $\Gamma(x)$, denoted by $\psi(x)=\frac{\Gamma'(x)}{\Gamma(x)}$, is called the psi or digamma function, and $\psi^{(k)}(x)$ for $k\in\mathbb{N}$ are called the polygamma functions. It is general knowledge that these functions are basic and that they have much extensive applications in mathematical sciences.
\par
In~\cite[Theorem~2]{Guo-Qi-TJM-03.tex}, the function
\begin{equation}\label{fun-ori}
\frac{[{\Gamma(x+y+1)}/{\Gamma(y+1)}]^{1/x}}{x+y+1}
\end{equation}
was proved to be decreasing with respect to $x\ge1$ for fixed $y\ge0$.
Consequently, the inequality
\begin{equation}\label{new3}
\frac{x+y+1}{x+y+2}\le\frac{[\Gamma(x+y+1)/\Gamma(y+1)]^{1/x}}
{[\Gamma(x+y+2)/\Gamma(y+1)]^{1/(x+1)}}
\end{equation}
holds for positive real numbers $x\ge1$ and $y\ge0$. Meanwhile, an open problem was posed in~\cite[p.~245]{Guo-Qi-TJM-03.tex} to ask for an upper bound $\sqrt{\frac{x+y}{x+y+1}}\,$ for the function in the right-hand side of the inequality~\eqref{new3}.
\par
In~\cite{Ya-Ming-Yu-JMAA-09}, the above-mentioned open problem was partially resolved as follows: If $y>0$ and $x>1$, then
\begin{equation}\label{yaming-ineq}
\frac{[\Gamma(x+y+1)/\Gamma(y+1)]^{1/x}}{[\Gamma(x+y+2)/\Gamma(y+1)]^{1/(x+1)}} <\biggl(\frac{x+y}{x+y+1}\biggr)^{1/2};
\end{equation}
if $y>0$ and $0<x<1$, then the inequality~\eqref{yaming-ineq} is reversed.
\par
For more information on the origin, history, backgrounds, motivations and recent developments of this topic, please refer to~\cite{Abram-Baric-Mat-Pecaric-JIA-07, alzer1, alzer2, alzer, Bennett-JIA-07, minc, Extension-TJM-2003.tex, Martin-JIA-07} and closely-related references cited therein.
\par
The aim of this paper is to extend the one-side inequality~\eqref{yaming-ineq} and to resolves the above-mentioned open problem.
\par
Our results may be stated as the following theorem.

\begin{thm}\label{Ya-Ming-Yu-extend-thm}
For $y+1>0$ and $x+y>0$, the inequality~\eqref{yaming-ineq} holds if $x>1$ and reverses if $x<1$. The cases $x=0,-1$ are understood to be the limits as $x\to0,-1$ on both sides of the inequality~\eqref{yaming-ineq}, that is,
\begin{equation}\label{yaming-ineq-no-x-1}
e^{\psi(y+1)}>(y+1)\biggl(\frac{y}{y+1}\biggr)^{1/2},\quad y>0
\end{equation}
and
\begin{equation}\label{yaming-ineq-no-x-2}
e^{-\psi(y+1)}>\frac1y\biggl(\frac{y-1}{y}\biggr)^{1/2},\quad y>1.
\end{equation}
Moreover, the powers $\frac12$ in~\eqref{yaming-ineq}, \eqref{yaming-ineq-no-x-1}, and~\eqref{yaming-ineq-no-x-2} are the best possible in the sense that the power $\frac12$ in the inequality~\eqref{yaming-ineq} can not be replaced by a larger number and that the powers $\frac12$ in the reversed inequality of~\eqref{yaming-ineq}, \eqref{yaming-ineq-no-x-1}, and~\eqref{yaming-ineq-no-x-2} can not be replaced by a smaller number.
\end{thm}

As a ready consequence of the proof of Theorem~\ref{Ya-Ming-Yu-extend-thm}, the following inequality is concluded.

\begin{cor}\label{Ineq-Ext-Cor}
For $x+y>0$ and $y+1>0$, if $0<\vert x\vert<1$, then
\begin{equation}\label{Ineq-Ext-Cor-ineq}
\biggl[\frac{\Gamma(x+y+1)}{\Gamma(y+1)}\biggr]^{1/x}>\biggr[\frac{(x+y)^{x+1}}{(x+y+1)^{x-1}}\biggr]^{1/2};
\end{equation}
if $\vert x\vert>1$, then the inequality~\eqref{Ineq-Ext-Cor-ineq} is reversed. In particular, the inequality
\begin{equation}\label{Ineq-Ext-Cor-ineq-y=0}
\Gamma(x+1)>\biggr[\frac{x^{x+1}}{(x+1)^{x-1}}\biggr]^{x/2}
\end{equation}
holds for $0<x<1$ and reverses for $x>1$.
\end{cor}

\section{Lemmas}

In order to prove Theorem~\ref{Ya-Ming-Yu-extend-thm}, we need the following lemmas.

\begin{lem}\label{wendel-gamma-ineq-lem}
For $k\in\mathbb{N}$ and $t>s>0$ with $t-s\ne1$, we have
\begin{equation}\label{wendel-gamma-ineq}
\min\biggl\{s,\frac{s+t-1}2\biggr\}<\biggl[\frac{\Gamma(s)}{\Gamma(t)}\biggr]^{1/(s-t)} <\max\biggl\{s,\frac{s+t-1}2\biggr\}
\end{equation}
and
\begin{equation}\label{n-s-ineq-st}
\frac{(k-1)!}{\bigl(\max\bigl\{s,\frac{s+t-1}2\bigr\}\bigr)^k}
<\frac{(-1)^{k-1} \bigl[\psi^{(k-1)}(t)-\psi^{(k-1)}(s)\bigr]}{t-s}
<\frac{(k-1)!}{\bigl(\min\bigl\{s,\frac{s+t-1}2\bigr\}\bigr)^k},
\end{equation}
where $\psi^{(0)}(x)$ stands for $\psi(x)$. Moreover, the lower and upper bounds in~\eqref{wendel-gamma-ineq} and~\eqref{n-s-ineq-st} are the best possible constants for which the inequalities hold.
\end{lem}

\begin{proof}
For real numbers $a$, $b$ and $c$, denote $\rho=\min\{a,b,c\}$, and let
\begin{equation}\label{h-def-sandor-new}
H_{a,b;c}(x)=(x+c)^{b-a}\frac{\Gamma(x+a)}{\Gamma(x+b)}
\end{equation}
with respect to $x\in(-\min\{a,b,c\},\infty)$. In~\cite[p.~283, Theorem~1]{sandor-gamma-3-note.tex-final}, it was obtained that
\begin{enumerate}
\item
the function $H_{a,b;c}(x)$ is logarithmically completely monotonic, that is,
$$
0\le(-1)^i[\ln H_{a,b;c}(x)]^{(i)}<\infty
$$
for $i\ge1$, on $(-\rho,\infty)$ if and only if
\begin{equation}\label{d1-dfn-new}
\begin{split}
(a,b;c)\in D_1(a,b;c)&\triangleq\{(a,b;c):(b-a)(1-a-b+2c)\ge0\}\\
&\quad\cap\{(a,b;c):(b-a) (\vert a-b\vert-a-b+2c)\ge0\}\\
&\quad\setminus\{(a,b;c):a=c+1=b+1\}\\
&\quad\setminus\{(a,b;c):b=c+1=a+1\};
\end{split}
\end{equation}
\item
so is the function $H_{b,a;c}(x)$ on $(-\rho,\infty)$ if and only if
\begin{equation}\label{d2-dfn-new}
\begin{split}
(a,b;c)\in D_2(a,b;c)&\triangleq\{(a,b;c):(b-a)(1-a-b+2c)\le0\}\\
&\quad\cap\{(a,b;c):(b-a) (\vert a-b\vert-a-b+2c)\le0\}\\
&\quad\setminus\{(a,b;c):b=c+1=a+1\}\\
&\quad\setminus\{(a,b;c):a=c+1=b+1\}.
\end{split}
\end{equation}
\end{enumerate}
See also \cite[pp.~1241\nobreakdash--1242, Theorem~4.1]{Feng-Qi-693.tex}. It is well-known that the limit
\begin{equation}\label{wendel-approx}
\lim_{x\to\infty}\biggl[x^{b-a}\frac{\Gamma(x+a)}{\Gamma(x+b)}\biggr]=1
\end{equation}
holds for real numbers $a$ and $b$, see~\cite[p.~257, 6.1.46]{abram} or~\cite[p.~3, Section~1.1.6]{bounds-two-gammas.tex}. This implies that
\begin{equation}\label{lim-H{a,b;c}(x)}
\lim_{x\to\infty}H_{a,b;c}(x)=1.
\end{equation}
\par
From the logarithmically complete monotonicity of $H_{a,b;c}(x)$, it is deduced that the function $H_{a,b;c}(x)$ is decreasing if $(a,b;c)\in D_1(a,b;c)$ and increasing if $(a,b;c)\in D_2(a,b;c)$ on $(-\rho,\infty)$. As a result of the limit~\eqref{lim-H{a,b;c}(x)} and the monotonicity of the function $H_{a,b;c}(x)$, it follows that the inequality $H_{a,b;c}(x)>1$ holds if $(a,b;c)\in D_1(a,b;c)$ and reverses if $(a,b;c)\in D_2(a,b;c)$, that is, the inequality
\begin{equation*}
x+\lambda<\biggl[\frac{\Gamma(x+a)}{\Gamma(x+b)}\biggr]^{1/(a-b)}<x+\mu
\end{equation*}
for $b>a$ holds if $\lambda\le\min\bigl\{a,\frac{a+b-1}2\bigr\}$ and $\mu\ge\max\bigl\{a,\frac{a+b-1}2\bigr\}$, which may be reduced to the inequality~\eqref{wendel-gamma-ineq} by replacing $x+a$ and $x+b$ by $s$ and $t$ respectively.
\par
Further, by virtue of the logarithmically complete monotonicity of $H_{a,b;c}(x)$ on $(-\rho,\infty)$ again and the fact in~\cite[p.~98]{Dubourdieu} that a completely monotonic function which is non-identically zero cannot vanish at any point on $(0,\infty)$, it is readily deduced that
\begin{multline*}
(-1)^k[\ln H_{a,b;c}(x)]^{(k)}=(-1)^k[(b-a)\ln(x+c)+\ln\Gamma(x+a)-\ln\Gamma(x+b)]^{(k)}\\
\begin{aligned}
&=(-1)^k\biggl[\frac{(-1)^{k-1}(k-1)!(b-a)}{(x+c)^{k}}+\psi^{(k-1)}(x+a)-\psi^{(k-1)}(x+b)\biggr]\\
&>0
\end{aligned}
\end{multline*}
for $k\in\mathbb{N}$ is valid if $(a,b;c)\in D_1(a,b;c)$ and reversed if $(a,b;c)\in D_2(a,b;c)$. Consequently, the double inequality
\begin{equation*}
-\frac{(k-1)!(b-a)}{(x+c_2)^{k}}<(-1)^k[\psi^{(k-1)}(x+b)-\psi^{(k-1)}(x+a)] <-\frac{(k-1)!(b-a)}{(x+c_1)^{k}}
\end{equation*}
holds with respect to $x\in(-\rho,\infty)$ if $(a,b;c_1)\in D_1(a,b;c)$ and $(a,b;c_2)\in D_2(a,b;c)$, which may be rearranged as
\begin{equation}\label{n-s-ineq}
\frac{(k-1)!}{(x+\alpha)^k}<\frac{(-1)^{k-1} \bigl[\psi^{(k-1)}(x+b)-\psi^{(k-1)}(x+a)\bigr]}{b-a} <\frac{(k-1)!}{(x+\beta)^k}
\end{equation}
for $x\in(-\rho,\infty)$ if $\alpha\ge\max\bigl\{a,\frac{a+b-1}2\bigr\}$ and $\beta\le\min\bigl\{a,\frac{a+b-1}2\bigr\}$, where $b>a$ and $k\in\mathbb{N}$. In the end, replacing $x+a$ and $x+b$ by $s$ and $t$ respectively in~\eqref{n-s-ineq} leads to~\eqref{n-s-ineq-st}. The proof of Lemma~\ref{wendel-gamma-ineq-lem} is thus complete.
\end{proof}

\begin{lem}\label{comp-thm-1}
For $x\in(0,\infty)$ and $k\in\mathbb{N}$, we have
\begin{equation}\label{qi-psi-ineq-1}
\ln x-\frac1x<\psi(x)<\ln x-\frac1{2x}
\end{equation}
and
\begin{equation}\label{qi-psi-ineq}
\frac{(k-1)!}{x^k}+\frac{k!}{2x^{k+1}}< (-1)^{k+1}\psi^{(k)}(x) <\frac{(k-1)!}{x^k}+\frac{k!}{x^{k+1}}.
\end{equation}
\end{lem}

\begin{proof}
In~\cite[Theorem~2.1]{Ismail-Muldoon-119} and \cite[Lemma~1.3]{sandor-gamma-2-ITSF.tex}, the function $\psi(x)-\ln x+\frac{\alpha}x$ was proved to be completely monotonic on $(0,\infty)$, i.e.,
\begin{equation}\label{com-psi-ineq-dfn}
(-1)^i\Bigl[\psi(x)-\ln x+\frac{\alpha}x\Bigr]^{(i)}\ge0
\end{equation}
for $i\ge0$, if and only if $\alpha\ge1$, and so is its negative, i.e., the inequality~\eqref{com-psi-ineq-dfn} is reversed, if and only if $\alpha\le\frac12$. In~\cite[Theorem~2]{chen-qi-log-jmaa}, \cite[Theorem~2.1]{Ismail-Lorch-Muldoon} and~\cite[Theorem~2.1]{Muldoon-78}, the function $\frac{e^x\Gamma(x)} {x^{x-\alpha}}$ was proved to be logarithmically completely monotonic on $(0,\infty)$, i.e.,
\begin{equation}\label{com-psi-ineq-ln-dfn}
(-1)^k\biggl[\ln\frac{e^x\Gamma(x)} {x^{x-\alpha}}\biggr]^{(k)}\ge0
\end{equation}
for $k\in\mathbb{N}$, if and only if $\alpha\ge1$, so is its reciprocal, i.e., the inequality~\eqref{com-psi-ineq-ln-dfn} is reversed, if and only if $\alpha\le\frac12$. Considering the fact in~\cite[p.~98]{Dubourdieu} that a completely monotonic function which is non-identically zero cannot vanish at any point on $(0,\infty)$ and rearranging either~\eqref{com-psi-ineq-dfn} or~\eqref{com-psi-ineq-ln-dfn} leads to the double inequalities~\eqref{qi-psi-ineq-1} and~\eqref{qi-psi-ineq}. Lemma~\ref{comp-thm-1} is proved.
\end{proof}

\begin{lem}[\cite{Mon-Two-Seq-AMEN.tex}]\label{Mon-Two-Seq-AMEN.tex-ineq}
If $t>0$, then
\begin{equation} \label{log-ineq-qi}
\frac{2t}{2+t}<\ln(1+t)<\frac{t(2+t)}{2(1+t)};
\end{equation}
If $-1<t<0$, the inequality~\eqref{log-ineq-qi} is reversed.
\end{lem}

\section{Proofs of Theorem~\ref{Ya-Ming-Yu-extend-thm} and Corollary~\ref{Ineq-Ext-Cor}}

Now we are in a position to prove Theorem~\ref{Ya-Ming-Yu-extend-thm} and Corollary~\ref{Ineq-Ext-Cor}.

\begin{proof}[Proof of Theorem~\ref{Ya-Ming-Yu-extend-thm}]
When $0\ge y>-1$ and $x>-y$, let
\begin{equation}\label{fy(x)ne0}
f_y(x)=\frac{\ln\Gamma(x+y+1)-\ln\Gamma(y+1)}x-\frac12\ln(x+y);
\end{equation}
When $y>0$ and $x>-y$, define
\begin{equation*}
f_y(0)=\psi(y+1)-\frac12\ln y
\end{equation*}
and $f_y(x)$ for $x\ne0$ to be the same one as in~\eqref{fy(x)ne0}.
Making use of the well-known recursion formula $\Gamma(x+1)=x\Gamma(x)$ and computing straightforwardly yields
\begin{multline}\label{fy(x+1)-fy(x)}
\begin{aligned}
f_y(x+1)-f_y(x)&=\biggl(\frac1{x+1}-\frac1x\biggr)\ln\frac{\Gamma(x+y+1)}{\Gamma(y+1)} \\ &\quad+\frac{\ln(x+y+1)}{x+1}+\frac12\ln\frac{x+y}{x+y+1}
\end{aligned}\\
=\frac1{x+1}\biggl\{\ln\biggl[\frac{(x+y)^{(x+1)/2}}{(x+y+1)^{(x-1)/2}}\biggr] -\ln\biggl[\frac{\Gamma(x+y+1)}{\Gamma(y+1)}\biggr]^{1/x}\biggr\}.
\end{multline}
\par
Substituting $s=y+1>0$ and $t=x+y+1>1$ into~\eqref{wendel-gamma-ineq} in Lemma~\ref{wendel-gamma-ineq-lem} leads to
\begin{equation*}
\min\biggl\{y+1,\frac{x+2y+1}2\biggr\} <\biggl[\frac{\Gamma(x+y+1)}{\Gamma(y+1)}\biggr]^{1/x} <\max\biggl\{y+1,\frac{x+2y+1}2\biggr\}
\end{equation*}
which is equivalent to
\begin{equation*}
\biggl[\frac{\Gamma(x+y+1)}{\Gamma(y+1)}\biggr]^{1/x} <
\begin{cases}
\dfrac{x+2y+1}2,& x>1\\
y+1,& x<1
\end{cases}
\end{equation*}
and
\begin{equation*}
\biggl[\frac{\Gamma(x+y+1)}{\Gamma(y+1)}\biggr]^{1/x} >
\begin{cases}
y+1,& x>1\\
\dfrac{x+2y+1}2,& x<1
\end{cases}
\end{equation*}
for $y+1>0$ and $x+y>0$. Consequently, it follows readily from~\eqref{fy(x+1)-fy(x)} that, for $y>-1$ and $x+y>0$,
\begin{enumerate}
\item
if $x>1$ and
\begin{equation}\label{x+2y+1}
\frac{(x+y)^{(x+1)/2}}{(x+y+1)^{(x-1)/2}}> \frac{x+2y+1}2,
\end{equation}
then $f_y(x+1)-f_y(x)>0$;
\item
if $-1<x<1$ and the inequality~\eqref{x+2y+1} reverses, then $f_y(x+1)-f_y(x)<0$.
\end{enumerate}
\par
For $x+y>0$ and $y>-1$, let
\begin{equation*}
g(x,y)=\frac{(x+y)^{x+1}}{(x+2y+1)^2(x+y+1)^{x-1}}.
\end{equation*}
The partial derivative of $g(x,y)$ with respect to $y$ is
$$
\frac{\partial g(x,y)}{\partial y}=\frac{1-x^2}{(x+2y+1)^3}\biggl(\frac{x+y}{x+y+1}\biggr)^x.
$$
This shows that
\begin{enumerate}
\item
when $\vert x\vert>1$, the function $g(x,y)$ is strictly decreasing with respect to $y>-1$;
\item
when $\vert x\vert<1$, the function $g(x,y)$ is strictly increasing with respect to $y>-1$.
\end{enumerate}
In addition, it is clear that $\lim_{y\to\infty}g(x,y)=\frac14$. As a result, it is easy to see that $g(x,y)\gtrless\frac14$ when $\vert x\vert\gtrless1$ for $x+y>0$ and $y>-1$. In other words, the inequality~\eqref{x+2y+1} is valid when $\vert x\vert>1$ and reversed when $\vert x\vert<1$ for all $x+y>0$ and $y>-1$. Consequently, the inequality $f_y(x+1)-f_y(x)>0$ holds if $x>1$ and reverses if $\vert x\vert<1$, where $x+y>0$ and $y>-1$.
\par
For $x<-1$, denote the function enclosed in the braces in \eqref{fy(x+1)-fy(x)} by $Q(x,y)$. Direct computation yields
\begin{multline*}
Q(x,y)=\frac{x+1}2\ln(x+y)-\frac{x-1}2\ln(x+y+1)-\frac1x\int_{y+1}^{x+y+1}\psi(u)\td u\\
=\frac{x+1}2\ln(x+y)-\frac{x-1}2\ln(x+y+1)-\int_0^1\psi((y+1)(1-u)+(x+y+1)u)\td u
\end{multline*}
and
\begin{align*}
\frac{\partial Q(x,y)}{\partial x}&=\frac{3 x+2 y+1}{2(x+y) (x+y+1)} +\frac12\ln\frac{x+y}{x+y+1}\\ &\quad-\int_0^1u\psi'((y+1)(1-u)+(x+y+1)u)\td u.
\end{align*}
Making use of the left-hand side inequality for $k=1$ in~\eqref{qi-psi-ineq} results in
\begin{multline*}
\frac{\partial Q(x,y)}{\partial x}<\frac{3x+2y+1}{2(x+y)(x+y+1)} +\frac12\ln\frac{x+y}{x+y+1}\\
-\int_0^1u\biggl\{\frac1{(y+1)(1-u)+(x+y+1)u}+\frac1{2[(y+1)(1-u)+(x+y+1)u]^2}\biggr\}\td u\\
=\frac12\biggl[\frac{x^2-2 y x-y (2 y+1)}{x (x+y) (x+y+1)}+\ln\frac{x+y}{x+y+1} -\frac{1+2y}{x^2}\ln\frac{y+1}{x+y+1}\biggr].
\end{multline*}
Further employing the left-hand side inequality of~\eqref{log-ineq-qi} in Lemma~\ref{Mon-Two-Seq-AMEN.tex-ineq} leads to
\begin{align*}
\frac{\partial Q(x,y)}{\partial x}&<\frac12\biggl[\frac{x^2-2 y x-y (2 y+1)}{x (x+y) (x+y+1)}-\frac2{1+2x+2y} +\frac{1+2y}{x^2}\cdot\frac{2x}{2+x+2y}\biggr]\\
&=\frac{(2y+3)x^2+2\bigl(y^2+2y+2\bigr)x+3y+2}{2(x+y)(x+y+1)(x+2y+2)(2x+2y+1)}\\
&\triangleq\frac{(2y+3)F_1(x,y)F_2(x,y)}{2(x+y)(x+y+1)(x+2y+2)(2x+2y+1)},
\end{align*}
where
\begin{equation*}
F_1(x,y)=\Biggl(x+\frac{2+y^2+2y-\sqrt{y^4+4 y^3+2 y^2-5 y-2}\,}{2y+3}\Biggr)
\end{equation*}
and
\begin{equation*}
F_2(x,y)=\Biggl(x+\frac{2+y^2+2y+\sqrt{y^4+4 y^3+2 y^2-5 y-2}\,}{2y+3}\Biggr).
\end{equation*}
For $x<-1$ and $x+y>0$, standard argument reveals that
\begin{equation*}
F_1(x,y)<\frac{2+y^2+2y-\sqrt{y^4+4 y^3+2 y^2-5 y-2}\,}{2y+3}-1<0
\end{equation*}
and
\begin{equation*}
F_2(x,y)>\frac{2+y^2+2y+\sqrt{y^4+4 y^3+2 y^2-5 y-2}\,}{2y+3}-y>0,
\end{equation*}
so $\frac{\partial Q(x,y)}{\partial x}<0$ and the function $Q(x,y)$ is decreasing with respect to $x<-1$. From the fact that $Q(-1,y)=0$, it follows that $Q(x,y)>0$ for $x<-1$, which means that when $y+1>0$ and $x+y>0$ the inequality~\eqref{yaming-ineq} is reversed for $x<-1$.
\par
For $x\gtrless1$, if
$$
\frac{[\Gamma(x+y+1)/\Gamma(y+1)]^{1/x}}{[\Gamma(x+y+2)/\Gamma(y+1)]^{1/(x+1)}} \gtrless\biggl(\frac{x+y}{x+y+1}\biggr)^{\alpha},
$$
then
$$
\alpha\lessgtr\frac{\ln\Gamma(x+y+1)-\ln\Gamma(y+1)}{x[\ln(x+y)-\ln(x+y+1)]} -\frac{\ln\Gamma(x+y+2)-\ln\Gamma(y+1)}{(x+1)[\ln(x+y)-\ln(x+y+1)]}
$$
is valid for $y+1>0$ and $x+y>0$. Since
\begin{align*}
&\quad\lim_{x\to1}\biggl\{\frac{\ln\Gamma(x+y+1)-\ln\Gamma(y+1)}{x[\ln(x+y)-\ln(x+y+1)]} -\frac{\ln\Gamma(x+y+2)-\ln\Gamma(y+1)}{(x+1)[\ln(x+y)-\ln(x+y+1)]}\biggr\}\\
&=\frac{\ln\Gamma(y+2)-\ln\Gamma(y+1)}{\ln(y+1)-\ln(y+2)} -\frac{\ln\Gamma(y+3)-\ln\Gamma(y+1)}{2[\ln(y+1)-\ln(y+2)]}\\
&=\frac12,
\end{align*}
it follows that $\alpha\lesseqgtr\frac12$. So the powers $\frac12$ in Theorem~\ref{Ya-Ming-Yu-extend-thm} are the best possible. Theorem~\ref{Ya-Ming-Yu-extend-thm} is thus proved.
\end{proof}

\begin{proof}[Proof of Corollary~\ref{Ineq-Ext-Cor}]
The inequality \eqref{Ineq-Ext-Cor-ineq} follows from the discussion in the proof of Theorem~\ref{Ya-Ming-Yu-extend-thm} about the positivity and negativity of the function enclosed by braces in~\eqref{fy(x+1)-fy(x)}.
\par
The inequality~\eqref{Ineq-Ext-Cor-ineq-y=0} is a special case of \eqref{Ineq-Ext-Cor-ineq} for $y=0$.
\end{proof}

\subsection*{Acknowledgemnts}
The authors are indebted to anonymous referees for their many valuable comments, corrections and suggestions.

\end{document}